\documentclass[12pt]
{article}

\usepackage{amsmath, amssymb}
\usepackage{amscd}
\usepackage[all]{xy}

\usepackage{enumerate}

\newcommand{\bZ}{\mathbb{Z}}
\newcommand{\cQ}{\mathcal{Q}}

\newcommand{\cE}{\mathcal{E}}
\newcommand{\cF}{\mathcal{F}}

\newcommand{\La}{\Lambda}

\begin{document}

\newtheorem{theorem}{Theorem}[section] % 1st argument is your name for it
\newtheorem{lem}[theorem]{Lemma}     
\newtheorem{cor}[theorem]{Corollary}
\newtheorem{prop}[theorem]{Proposition}
\newtheorem{remark}[theorem]{Remark}
\newtheorem{defi}[theorem]{Definition}
\newtheorem{expl}[theorem]{Example}
\newtheorem{rem}[theorem]{Remark}
\newtheorem{notation}[theorem]{Notation}

\parindent0pt

\baselineskip=16pt

\begin{center}{\bf\large Algebras with non-periodic bounded modules} \\
\medskip
\textit{Karin Erdmann}
\medskip
\end{center}

\begin{abstract} We study weakly symmetric special biserial algebra
of infinite representation type.
We show that usually some socle deformation of 
such an algebra  has non-periodic bounded modules.  The exceptions are
precisely the algebras whose 
Brauer graph is a tree
with no multiple edges.
If the algebra has a non-periodic bounded module then its 
 Hochschild cohomology    cannot satisfy 
the finite generation property (Fg) introduced in \cite{EHSST}. 
\end{abstract}

AMS Subject classification: 16G20, 16E40.

\medskip

\section{Introduction}

\ \ Assume $\La$ is a finite-dimensional selfinjective algebra over some field $K$. If $M$ is a finite-dimensional non-projective $\La$-module, let 
$$\ldots \to P_n\stackrel{d_n}\to P_{n-1} \to \ldots \to P_1\stackrel{d_1}\to P_0
\stackrel{d_0}\to M\to 0
$$
be a minimal projective resolution of $M$.  The 
module $M$ is called bounded if the dimensions of the projectives $P_n$
have a common upper bound,  that is,  $M$ has complexity one. 
The kernel of $d_n$ is  the syzygy  $\Omega^n(M)$;   we say that the module
$M$ is periodic if $\Omega^d(M)\cong M$ for some $d\geq 1$.
A periodic module has complexity one but the converse need not hold. 
We call a module   $M$  a {\it criminal }  if it has complexity one but is not periodic.
We would like to understand which algebras have criminals.

\medskip

\ \ J. Alperin proved in \cite{Alp} that the 
group algebra of a finite group does not have criminals  when the coefficient field is 
algebraic over its prime field. 
On the other hand, R. Schulz discovered that there are four-dimensional selfinjective algebras 
which have criminals, see \cite{Sch}. 
In the context of commutative algebra,  there is a similar problem.
Eisenbud proved in \cite{Ei} that for  complete intersections,  if  a finitely generated
module has bounded Betti numbers  then it is eventually periodic. He conjectured that this should be true for any
commutative Noetherian local ring. However, counterexamples
were constructed  by Gasharov and Peeva \cite{GP}.

 \medskip
 
 \ \ Subsequently, a theory of support varieties was developed for modules of group algebras of
 finite groups. This is based on group cohomology and depends crucially on that the fact that it is Noetherian.
 It follows from this theory that a group algebra over an arbitrary field does not have criminals, so
 that Alperin's theorem holds in general; for a proof see 2.24.4 in \cite{Be1}.
 More recently, a support variety theory was developed  for modules of selfinjective algebras,  based on Hochschild cohomology \cite{SnS, EHSST}. 
 This also requires suitable finite generation, namely
 the Hochschild cohomology $HH^*(\La)$ should be Noetherian and the ext-algebra  of $\La$ should be 
finitely generated as a module over $HH^*(\La)$. This condition is called (Fg) in \cite{So}, it is equivalent to (Fg1, 2) in
\cite{EHSST}. 
Again, if (Fg) holds for $\La$, so that  $\La$-modules have support varieties, then $\La$ does not have criminals 
(see 5.3 in \cite{EHSST}).

\bigskip

\ \ 
The algebras studied by Schulz therefore do not satisfy (Fg).   More generally, 
weakly symmetric algebra with radical cube zero were investigated in \cite{ES}, \cite{E2}.  
The algebras in these papers   which have criminals happen to be   special biserial,  therefore one may ask  when a special biserial
weakly symmetric algebra has criminals. Of course, if an algebra has a chance to have criminals
it must have infinite representation type.

\bigskip

\ \ Here we study special biserial weakly symmetric $K$-algebras of infinite representation type, we assume  $K$ is an algebraically closed
field which contains non-roots of unity. 
An algebra is special biserial weakly symmetric if its basic algebra satisfies 2.1. Existence of criminals is invariant under Morita equivalence,
and  we will work throughout with basic algebras. We  assume that the algebra is indecomposable, so that its quiver is connected.
The algebras in 2.1 have socle relations involving scalar parameters, so that we have a family of algebras,  which we write as $\La_{\bf q}$
where ${\bf q}$ is the collection of the the socle scalars. Each algebra in a family is a socle deformation of the algebra $\La_{\bf 1}$ for which all socle
scalars are equal to $1$. 

\ \ 
Recall that if  $\La$ and $\Gamma$ are selfinjective, then $\Gamma$ is a socle deformation of $\La$ if $\Gamma/{\rm soc}(\Gamma)$ is isomorphic to
$\La/{\rm soc}(\La)$.
For example when the field has characteristic $2$ then the algebra studied by Schulz is a socle deformation of
the group algebra of a Klein 4-group.  There are similar socle deformations for group algebras of dihedral 2-groups, which 
are also special biserial and weakly symmetric.

\medskip
\ \  Our main result  answers   when there is some choice for ${\bf q}$ such that the algebra $\La_{\bf q}$ has criminals. 
The algebra $\La_{\bf q}$ has a Brauer graph $G_{\La}$, which is is independent of ${\bf q}$,  we  define it in section 4. 
The Brauer graph generalises the Brauer tree for a block of a group algebra of finite type.  For the algebras of infinite type, 
this graph is usually not a tree.  We will prove the following:

\bigskip
\begin{theorem} \label{maintheorem} Let $K$ be an algebraically closed
field which contains some non-roots of unity. 
Assume $\La_{\bf q}$ is a family of 
indecomposable weakly symmetric 
special biserial $K$-algebras of infinite type. Then the following are equivalent:\\
(a)  $\La_{\bf q}$ does not have criminals for any ${\bf q}$.\\
(b) The Brauer graph $G_{\La}$  is a tree with no multiple edges.
\end{theorem}

\bigskip

From the perspective of group representation theory, one might have expected that  \lq most \rq \   selfinjective  algebras should satisfy (Fg) and have a support variety theory.
Our Theorem suggests  that this may  not be the case.

\ \ 
A special biserial weakly symmetric algebra $\La_{\bf q}$ of finite type does not  have criminals since it cannot have
a non-periodic module, and its  Brauer graph is a tree. As well, such an  algebra is symmetric and is isomorphic to $\La_{\bf 1}$.
For such an algebra, there is a unique $\Omega$-orbit consisting of maximal uniserial
modules and simple modules whose projective cover is uniserial. In the case of a block of a group algebra, this is due to J.A. Green \cite{G}, 
but it holds for arbitrary symmetric indecomposable Brauer tree algebras of finite type.
Assume   $\La_{\bf q}$ 
is weakly symmetric of infinite type whose  Brauer graph is a tree. Then
there is also  a unique $\Omega$-orbit consisting of maximal uniserial modules  generated by arrows and simple modules whose
projective cover is uniserial,
this  follows from \cite{Ro} and Section 4.1 in \cite{DD} when the algebra is symmetric. 

\bigskip

\ \  Given the presentation of the algebra as in 2.1, it is easy to determine its Brauer graph.  Our Theorem shows 
that the group algebra of a dihedral 2-group over characteristic 2 has a socle deformation with criminals, as well the group algebra
of the alternating group $A_4$, or also other special biserial algebra occuring in algebraic Lie theory, see for example 
\cite{FSk}.  Algebras of dihedral type, as defined in \cite{E1}, are special biserial and symmetric, hence they are examples for the algebras
in the Theorem. Of the ones which occur as blocks,  only one has a Brauer graph which is a tree. The Hecke algebras of type $A$ which have tame 
representation type are also special biserial, there are two Morita
equivalence classes, see \cite{EN}. Both have Brauer graphs which are trees.

\bigskip

\ \ 
In Section 2, we define the algebras and summarise  properties we need. In Section 3 we recall the definition of band modules, and
we determine the $\Omega$-translates of certain band modules. 
Using this we prove the Theorem in Section 4. 

\ \ 
It would  be interesting to know whether 
the algebras in Theorem 1.1 for which the Brauer graph is a tree, always satisfy (Fg). 
More generally, one may ask whether any selfinjective algebra
where modules have finite complexity, and which does not have criminals,
must satisfy (Fg).

\newpage

\section{The algebras}

Assume  $\La = K\cQ/I$ where $\cQ$ is a finite connected 
quiver and $I$ is an admissible ideal
of the path algebra $K\cQ$, generated by a set of relations $\rho$.

\begin{defi} \normalfont \label{Def:SB} 
The algebra $\La$ is special biserial and weakly symmetric
 if it satisfies the following: 
\begin{enumerate}
\item 
Any vertex of $\cQ$ is either the source of  two arrows and
is the target of two arrows, or it is the source of one arrow and is
the target of one arrow. We say that the vertex has valency two, or one, respectively.
\item 
If two  different arrows $\alpha$ and $\beta$ start at vertex $i$, then for
 an arrow $\gamma$ ending at $i$, precisely one of the two paths 
$\gamma\alpha$ or $\gamma\beta$ is in $\rho$.
\item 
If two different arrows $\gamma$ and $\delta$ end at vertex $i$, then
for  an arrow $\alpha$ starting at $i$, precisely  one of the two paths
$\alpha\gamma$ or $\alpha\delta$ is in $\rho$.
\item
For each vertex $i$ of $\cQ$ of
valency two, there are different paths $C_i$ and $D_i$ of length $\geq 2$
starting and ending 
at $i$, and
non-zero scalars $p_i, q_i \in K$ such that 
$p_iC_i + q_iD_i$ belongs to $\rho$. 
\item 
For each vertex $i$ of $\cQ$ of valency one, there is a path $C_i$ of length
$\geq 2$ such that $C_i\alpha$ belongs to $\rho$ where $\alpha$ is the
arrow starting at $i$. 
\item Any rotation of a path $C_i$ or $D_i$ in 4. or 5. is a path occuring
in the relation of the form 4. or 5. 
\item The set $\rho$ consists precisely of the relations described above.
\end{enumerate}
The algebra depends on the parameters $p_i, q_i$, and we write 
$\La = \La_{\bf q}.$
\end{defi}

We refer to relations  4. and 5. as 'socle relations'.
The definition of a special biserial algebra is slightly more general,
details may be found for example  in \cite{BR}.  
%We assume for the next few sections that all vertices have valency two.
%This means that the relations 5. do not occur.
%Later  we will explain how this also deals with
%the general case. 
If such an algebra has infinite representation type then there must
be at least one vertex with valency two.

\begin{rem} \normalfont 

(1) \ We identify as usual paths with the image in the algebra. So for example in $\La$ we have
$p_iC_i + q_iD_i=0$ but $C_i$ and $D_i$ are non-zero. 
The element $C_i$ spans the  
 socle of the indecomposable projective $e_i\La$, and using this, one can write down a non-degenerate
 bilinear form on $\La$ verifying that the algebra is indeed  selfinjective. Hence it is weakly symmetric, noting that the simple quotient and the socle of 
 $e_i\La$ are isomorphic.

(2) Let $e_i$ be a vertex of $\cQ$. Then the indecomposable
projective module  $e_i\La$ 
has a  basis consisting of  all proper initial subwords of
$C_i$ and $D_i$ of positive length, together with $e_i$ and one of $C_i$ or
$D_i$. In particular $\dim e_i\La = |C_i| + |D_i|$ where we write $|\eta|$ for
the number of arrows in  $\eta$.  
 \end{rem}

\bigskip

\begin{defi}\label{sigma}  Given a special biserial weakly symmetric
algebra, there is an associated permutation $\sigma$ of the
arrows of $\cQ$: \  
For each arrow $\alpha$ of $\cQ$, define $\sigma(\alpha)$ to 
be the unique arrow such that $\alpha\dot\sigma(\alpha)$ is non-zero in$\La$. 
\end{defi}
We may write $\sigma$  as a
product of disjoint cycles. Any monomial $C_i$ occuring in
a socle relation is then the product
$(\alpha_1\alpha_2\ldots \alpha_r)^m$ where $(\alpha_1 \ \alpha_2 \ \ldots \ \alpha_r)$ is a cycle of $\sigma$, and $m\geq 1$. 
Similarly $D_i= (\beta_1\beta_2\ldots \beta_s)^t$ where the product of the
$\beta_j$ is taken over a cycle of $\sigma$. This cycle may or may not 
be the same as the cycle
$(\alpha_1 \ \alpha_2 \ \ldots \ \alpha_r)$.

\subsection{Candidates for criminals}

The indecomposable non-projective $\La$-modules are classified, they are \lq strings \rq \ or \lq bands \rq . A description, and further details, may be
found in \cite{BR} or \cite{E1}. A candidate to be a criminal must have complexity one.

(1) If $M$ is a string module and is not of the form $\alpha \La$ or $e\La/\alpha\La$ for an arrow $\alpha$ starting at vertex $e$, or is in the Auslander-Reiten component of such module,  then
$M$ has complexity $\geq 2$. One can see this for example by considering its Auslander-Reiten translates $\tau^rM$ for $r\geq 2$. For a 
selfinjective algebra, the Auslander-Reiten translation $\tau$ is isomorphic to $\Omega^2\circ \nu$ where $\nu$ is the
Nakayama automorphism of the algebra, see for example \cite{ARS}, IV.3.7. As one can see from the construction of
irreducible maps, the dimensions of $\tau^r(M)$ are unbounded for $r\geq 1$. 
Hence  the dimensions of the modules $\Omega^{2r}(M)$ are also
unbounded, which implies that
$M$ has complexity $\geq 2$. 

(2) There are finitely many Auslander-Reiten components containing string modules of the form $\alpha\La$ or $e\La/\alpha \La$. 
These string modules 
are permuted by   $\Omega$ and hence are 
$\Omega$-periodic, since the set of arrows is finite. 
The action of $\Omega$ induces an equivalence of the stable module category,
and it commutes with $\tau$,   and
it follows that all modules in these components must be periodic with respect
to $\Omega$, and cannot be  criminals.

\medskip

(3) Band modules have  complexity one. They are parametrized by a band word $W$, a non-zero scalar $\lambda$, and 
a non-zero vector space $V$, we give details below. 
If $\lambda$ and $W$ are fixed, the corresponding band modules as $V$ varies, form one Auslander-Reiten component. 
Again, since  $\Omega$ induces an equivalence of the stable module category and commutes with $\tau$, the component contains a criminal
if and only if the band module with $V=K$ is a criminal. Therefore we can focus on band modules 
where the space is $K$.

(4) We also note that a special biserial algebra
has infinite representation type if and only if there are band modules.
This is proved in
\cite{SW}, Theorem 1 and Lemma 2 (with the terminology of primitive
V-sequences, for translating terminology, see \cite{WW}).

\section{Band modules}

We are looking for criminals, and therefore we focus
on 
band modules. We  start by describing the parameter set.
It is convenient to identify a vertex $i$ of $\cQ$ with the corresponding
idempotent $e_i$ of the path algebra $K\cQ$. 

\bigskip

\begin{defi}\label{bandword}\normalfont Let $e_0, \ldots, e_m$ and
$f_0, \ldots, f_m$ be vertices in $\cQ$ of valency two, and let $e_{m+1}=e_0$.
 A band word $W$ is a sequence
$(a_i, b_i)_{i=0}^m$ where the $a_i$ and $b_i$ are paths in $\cQ$ between 
vertices of valency two, where $a_i: e_i \mapsto f_i$ for $0\leq i\leq m$, 
and $b_i: e_i\mapsto f_{i-1}$ for $1\leq i\leq m+1$, such that
$a_i$ and $b_{i-1}$ are proper initial subpaths of the $C_i$ and $D_i$. 
Moreover, the sequence $(a_i, b_i)$ must be minimal with these properties.

That is, there is no shorter sequence $(\tilde{a}_i, \tilde{b}_i)$ with
the same properties such that $(a_i, b_i)$ is the concatenation 
of copies of $(\tilde{a}_i, \tilde{b}_i)$.

The band word $W$ may be described by a 
quiver:
$$e_0 \stackrel{a_0} \longrightarrow  f_0 \stackrel{b_0} \longleftarrow e_1 \stackrel{a_1}\longrightarrow f_1 \longleftarrow \ldots f_m \stackrel{b_m}\longleftarrow e_{m+1}= e_0
$$
Note that we do not specify the 
 the names of the arrows occuring in the paths $a_i, b_i$, since we will not
need these.
For details, we refer to \cite{BR} or \cite{E1}.
\end{defi}

For example, if  all vertices of the quiver have valency two, then
there is such band word where all the $a_i$ and $b_i$ are arrows. 
In this case, the minimality condition holds precisely if all the $e_i$ are distinct, equivalently if all the $f_i$ are distinct.

\begin{defi}\label{bandmodule}
%\normalfont 
The band  module $M(\lambda)$
 associated to the band word
$W$ as in \ref{bandword} and a vector space $V$, labelled
by a parameter $0\neq \lambda \in K$, is 
defined as follows: 
\begin{enumerate}
\item
For each vertex along the paths $a_i$ and $b_i$,  except the for the start
vertex of $b_m$, 
we take a copy of $V$. We identify the space at the start of $b_m$ with the
space at the start of $a_0$.  
\item The first arrow of $a_0$ acts by multiplication with an indecomposable
Jordan block matrix with eigenvalue $\lambda$.
\item All other arrows  occuring 
in the paths $a_i, b_i$ act as identity. 
\end{enumerate}
\end{defi}

We will only take $V=K$, then the first arrow of $a_0$ is multiplication
by $\lambda$. The 
module has dimension $\sum_{i=0}^m |a_i| + |b_i|$ where
$|\eta|$ is the number of arrows in the path $\eta$. It is 
 indecomposable, and
$M(\lambda)\cong M(\mu)$ only if $\lambda =\mu$. 

\begin{rem} \normalfont The arrow which acts by a non-identity scalar need not be
the first arrow of $a_0$. There are variations which give isomorphic
modules, details are discussed in  
\cite{BR}, or \cite{WW}. 
\end{rem}

\bigskip

\begin{expl} \normalfont \
To illustrate  the shorthand notation, let $m=0$ and
$a_0 = \alpha_1\alpha_2\alpha_3$, and $b_0=\beta$, then 
the word written in detail is
$$e_0 \stackrel{\alpha_1}\rightarrow \cdot \stackrel{\alpha_2}\rightarrow
 \cdot \stackrel{\alpha_3} \rightarrow  f_0 \stackrel{\beta}\leftarrow e_0.
$$
The   module $M(\lambda)$ as defined in \ref{bandmodule}, associated to  this word and $V=K$, is four-dimensional.
\end{expl}

\bigskip

We fix the word $W$, and the module $M(\lambda)$, and
we will now determine $\Omega$-translates for $M(\lambda)$. Note that
$\Omega^2(M(\lambda))$ will be a band module defined by the same word $W$, and
therefore we only need to calculate two steps.
This requires 
using the socle relations for the vertices $f_t$ and $e_t$ occuring
in the word $W$. We fix now the notation for these, so that we can keep track
over the
paths $a_t$ and $b_t$.

\bigskip

\begin{notation}\label{notation} \normalfont 
We write 
the socle relation relation starting and ending
at the vertex $f_t$ in the form 
$$ (\theta_t) \ \ \ \ \ \ p_t(A_ta_t) + q_t(B_tb_t) = 0, 
$$
where $A_t$, $B_t$ are paths, and $p_t, q_t$ are non-zero scalars.
Here 
$$A_t: \  f_t \longrightarrow \ e_t, \ \ B_{t}: f_{t} \longrightarrow \ e_{t-1}
$$
(taking indices modulo $m+1$). 
Similarly we write the socle relations starting and ending at vertex $e_t$
in the form
$$(\theta_t') \ \ \ \ \ \ p_t'(a_tA_t) + q_t'(b_{t-1}B_{t-1})=0,
$$
where $p_t'$ and  $q_t'$ are non-zero scalars.
\end{notation}

\begin{rem} \normalfont

Note that with this notation, we have 
$a_tB_t=0= b_tA_t$ for each $t$. For example,
let $\alpha$ be the last arrow of $a_t$. Since $a_tA_t$ is non-zero,
we know that $A_t$ must start with $\sigma(\alpha)$ where $\sigma$ is the
permutation defined in \ref{sigma}. The first arrow
of $B_t$ is the other arrow starting at $f_t$, say this is $\beta$,  and by condition 2. of
Definition 2.1 we have $\alpha \beta=0$ and hence $a_tB_t=0$. 
Similarly  
we have $A_tb_{t-1}=0 = B_ta_{t+1}$. 
\end{rem}

\bigskip

\begin{prop}\label{Omega} Let $v:= \prod_{t=0}^{m} (q_i/p_1) \prod_{t=0}^m (p_i'/q_i')$.
Then $\Omega^2(M(\lambda))$ is isomorphic to $ M(v\lambda)$. 
\end{prop}

\bigskip
\begin{rem} \normalfont 
Hence we have 
that 
$\Omega^{2r}(M(\lambda)) \cong M(v^r\lambda)$ for $r\geq 1$. 
This shows directly that $M(\lambda)$ has a bounded projective resolution, that
is, its complexity is one.
If $v=1$ then $M(\lambda)$ has $\Omega$-period at most two, 
and this occurs when the algebra is symmetric. 
As well $v$ might be some root of unity but $v\neq 1$. If so then $M(\lambda)$ is still periodic but it can have a larger period. Then we see that
$\Omega^2(M(\lambda))$ is not isomorphic to $\tau(M(\lambda))$ and we deduce that the Nakayama automorphism is non-trivial.
Our main interest here is in algebras for which $v$ is not a root of unity.
Note that the parameter $v$ depends only on the band word $W$ but not on $\lambda$. We  say that '$v$ is the parameter for $W$'.
\end{rem}

\subsection{The case $m=0$}

We prove the Proposition first for a band word with $m=0$, this 
needs slightly different (and less) notation. 
In this case we have 
 paths $a, b: e\mapsto f$ and the socle relations at $e$ and $f$ 
are of the form 
$$p(Aa) + q(Bb) = 0, \ \ p'(aA) + q'(bB)=0
$$
where $p, q, p'$ and $q'$ are non-zero scalars.

Fix some non-zero $\lambda\in K$, we want to 
construct the band module $M(\lambda)$ 
as a submodule of $f\La$.  That is, we look for an element
$w\in f\La$ such that $wa = \lambda\cdot wb$ and such that
this is a non-zero element in the socle of $f\La$, ie, it is a non-zero
scalar multiple of $Aa$. 

\bigskip

\begin{defi} \normalfont   Let 
$w:= \lambda A  - \frac
{q}{p}B \in f\La.$
Since $Ba=0$ and $Ab=0$ (see \ref{notation}), we have
$$wa = \lambda Aa = -\lambda(q/p)Bb \mbox{ and }
wb = -(q/p)Bb
$$ and hence $wa=\lambda wb$, and this is non-zero in the socle.
Let $a= \alpha_1\alpha_2\ldots \alpha_r$ where
the $\alpha_i$ are arrows, and let $b=\beta_1\ldots \beta_s$ for
arrows $\beta_j$. We may write down
a basis for  $w\La$, where each basis vector spans the
1-dimensional space at a vertex of the quiver described
in \ref{bandword},  showing that $w\La$ is of the form as in 
\ref{bandmodule}.
Namely, take the basis
$$w, \ \lambda^{-1}w\alpha_1, \ \lambda^{-1}w\alpha_1\alpha_2, \ldots, 
\lambda^{-1}wa, \ \ w\beta_1, \ w\beta_1\beta_2, \  w\beta_1\beta_2\ldots
\beta_{s-1}.
$$
Hence $w\La$ is isomorphic to $M(\lambda)$.
\end{defi}

\bigskip
\subsubsection{ The module $\Omega(M(\lambda))$.}

We find $\Omega(M(\lambda))$, this can be identified with the kernel of
the homomorphism
$$\psi: e\La \to w\La, \ \ \psi(x) := wx.
$$
We see that
$$w(a - \lambda b) = \lambda Aa + (q/p)\lambda Bb  = 0.
$$
Hence if $\zeta:= a - \lambda b \in e\La$, then 
$\zeta \La$ is a submodule of $\Omega(M(\lambda))$. 
We compare  dimensions; the dimension of  
$\zeta \La$ is  $|A| + |B|$. As well the dimension of $w\La$ is $|a| + |b|$ and hence
the sum of the dimension is equal to the dimension of $e\La$.
It follows that
$\zeta \La = \Omega(M(\lambda))$. 

\medskip

\subsubsection{The module $\Omega^2(M(\lambda))$.}

We identify $\Omega^2(M(\lambda)) \cong \Omega(\zeta\La)$
with the kernel of the map $\psi^+: f\La \to \zeta\La$ given by 
left multiplication with $\zeta$. 
Let $w^+:= \lambda A - (q'/p')B$, then 
$$(a -\lambda b)w^+ = \lambda aA + \lambda (q'/p')bB = 0.
$$
As before, we compare dimensions and deduce that $w^+\La = \Omega^2(M(\lambda))$. 

We identify $w^+\La$. 
First we have $w^+a = \lambda Aa$ and $w^+b = -(q'/p')Bb= Aa$ and hence
$$w^+a = \lambda v(w^+b).$$
As well, this is a non-zero element in the socle of $f\La$. Hence
$\Omega^2(M(\lambda) \cong M(v\lambda)$ where
$v = (q/p)(p'/q')$, as stated in the Proposition.

\bigskip

\begin{expl} \label{example:m=0}\normalfont

(1)  Let $\La$ be the local algebra with generators $x, y$ and relations
$$x^2=0=y^2, \  \ p(yx)^2 + q(xy)^2 = 0
$$
where $p, q$ are non-zero scalars. 

We have the band word $W$ given by $a=x$ and $b=y$. 
The relevant 
socle relations are then
$$p(Ax) + q(By)=0 \ \mbox{ and } \  \ p'(xA) + q'(yB)=0
$$
with $A=yxy$ and $B=xyx$ and $p'=q, q'=p$. Therefore we have that
$v=(q/p)^2$. 
If $q/p$ is not a root of unity then the modules $M(\lambda)$ are
criminals for the algebra with these parameters.

If char$(K)=2$, then the algebra with $q=p=1$ is isomorphic to the group 
algebra of the dihedral group of order $8$. 

\bigskip

(2) \ 
There is  family of commutative special biserial local algebras with 
generators $x, y$, and relations
$$xy=0 = yx, \  \ p(x^r) + q(y^s) =0
$$
for $r, s\geq 2$ and $p, q$ non-zero scalars.
We have again the band word $W$ given by $a=x$ and $b=y$. 
Writing down the socle relations in this case, we see that the parameter
 $v$ is equal to $1$ in this case. 

%We see that we may re-scale the generators. Set $x'= cx$ and $y'=dy$ where
%$c$ is a root of $t^r-p^{-1}$ and $d$ is a root of $t^s-q^{-1}$. Then 
%the elements $x'$ and $y'$ also are generators of the algebra, and 
%the algebra is given by the relations 
%$$x'y'=0, y'x'=0, \ (x')^r + (y')^s = 0.
%$$
%This shows that any algebra $\La_q$ is isomorphic to the algebra $\La_1$
%with all parameters equal to $1$. 

\medskip

(3)  Let $\La$ be the algebra with quiver

\[
\xymatrix{
0\ar@(ul,dl)_{\alpha} \ar@<1ex>[r]^-{\beta} &
1\ar@<1ex>[l]^-{\gamma}
}
\]
and relations
$$p\alpha^2 + q(\beta\gamma)^s=0, \ \ \alpha\beta=0, \gamma\alpha=0.
$$
where $p, q\neq 0$.  
Take the words  $a = \alpha: e_0 \rightarrow e_1$ and $b = \beta\gamma:
e_0\rightarrow e_1$, then $(a, b)$ is a band word. In this case $A=\alpha$ and $B = (\beta\gamma)^{s-1}$. The two socle relations we need in this case
are identical since $aA = Aa$ and $bB= Bb$. The parameter $v$ is equal to $1$.
When $s=2$, this algebra  occurs as a tame Hecke algebras, see
\cite{EN}. 
\end{expl}

\bigskip

\subsection{The case $m\geq 1$}

We take a band word $W$ as described in \ref{bandword}, with $m\geq 1$,  and 
we will construct $M(\lambda)$ by specifying generators, 
as a submodule of $\oplus_{t=0}^m f_t\La$. 

\newpage

\begin{defi} Define elements in the direct sum $\oplus_{t=0}^m
f_t\La$:
\begin{align*} v_0:= & (c_0A_0, 0, \ldots \ ,d_mB_m)  \cr 
v_1:= &(d_0B_0, c_1A_1, 0, \ldots , 0)& \ \cr
v_2:=& (0, d_1B_1, c_2A_2, 0, \ldots \ )& \ \cr
& \ldots & \ \cr
v_m:= & (0, \ldots, 0, d_{m-1}B_{m-1}, c_mA_m)&\
\end{align*}
where the  $c_i$ and the $d_i$ are  non-zero scalars.
\end{defi}

With this, we have
\begin{align*}v_0a_0=& (c_0(A_0a_0),  0, \ldots, 0) \cr
v_1b_0=&(d_0(B_0b_0), 0,  \ldots, 0)\cr
v_1a_1 = &(0, c_1(A_1a_1), 0,   \ldots)\cr
v_2b_1=& (0, d_1(B_1b_1), 0,  \ldots )
\cr
&\ldots\cr
v_ma_m=&(0, \ldots, 0, c_m(A_ma_m))\cr
v_0b_m=&(0, 0, \ldots, d_m(B_mb_m)).
\end{align*}

For any choice of scalars $c_i$ and $d_i$, the elements $v_i$ generate a submodule of $\oplus_{i=0}^m f_t\La$, and its dimension depends only on
the length of the $a_i$ and $b_i$ and one finds that
the dimension is equal to $\dim M(\lambda)$. 
We can see from the parameters when it is isomorphic to $M(\lambda)$.

\begin{lem}\label{lem:ident}
The submodule of $\oplus_{t=0}^m f_t\La$ generated by 
$v_0, v_1, \ldots, v_m$ is isomorphic to $M(\lambda)$ if and only if
$$\lambda p_0d_0 + q_0c_0 = 0  \  \  \mbox{ and } \ \  p_td_t+ q_tc_t=0
$$
for $1\leq t\leq m$.
\end{lem}

\bigskip

{\it Proof } We need $v_0a_0=\lambda v_1b_0$, that is
$$c_0(A_0a_0) = \lambda d_0(B_0b_0)
$$
We have $(B_0b_0) =-(p_0/q_0)(A_0a_0)$ and $A_0a_0\neq 0$. 
Substituting this gives the first equation. 
Similarly we need $v_1a_1=v_0b_1$, that is
$$c_1(A_1a_1)  = d_1(B_1b_1)
$$
Using $B_1b_1= -(p_1/q_1)(A_1a_1)$ gives the second equation.
Similarly the other equations follow.
Conversely, if all these identities hold then 
$\sum_{t=0}^m v_t\La \cong  M(\lambda)$. 
$\Box$

\vspace*{2cm}

We continue with the notation
$$c_0=\lambda p_0, \ c_i = p_i, \ \ d_j=-q_j \ (1\leq i\leq m, 0\leq j\leq m).
$$
We construct the first two steps of a minimal projective resolution.

\bigskip

\subsubsection{The module $\Omega(M(\lambda))$.}

Let $\Psi: P_0= \oplus_{t=0}^m e_t\La
\longrightarrow \oplus_{t=0}^m f_t\La$ 
be the
map given by left multiplication with the matrix 
$$\left(\begin{matrix} \lambda p_0A_0 & -q_0B_0 & 0 &  \ldots&0 & 0\cr
0  &   p_1A_1 & -q_1B_1 &  0&\ldots &0\cr
\ldots & & & &&\cr
0&0&\ldots&0 &p_{m-1}A_{m-1}& -q_{m-1}B_{m-1}\cr
-q_mB_m & 0 & \ldots&0 &0& p_mA_m\end{matrix}\right).
$$
Then $\Psi$ takes the standard generators of $P_0$ to
$v_0, v_2, \ldots, v_m$, and hence the image of $\Psi$ is
$M(\lambda)$. 
We know that $\Omega (M(\lambda))$ has minimal projective cover of the form
$$P_1= \oplus_{t=0}^m f_t\La \stackrel{\Psi_1}\to P_0$$
and 
$\Psi_1$ is given by left multiplication with a matrix of the form
$$\left(\begin{matrix} r_0a_0 & 0& 0& \ldots & 0& s_mb_m\cr
                       s_0b_0 & r_1a_1& 0&\ldots &0&0\cr
                         0&     s_1b_1& r_2a_2&\ldots &0&0\cr
                        \ldots &&&&&\cr
                       0&0&0&\ldots &s_{m-1}b_{m-1} & r_ma_m\end{matrix}
\right). 
$$
Here the $r_t$ and $s_t$ are non-zero scalars. By comparing dimensions, we
see that ${\rm Im}(\Psi_1) = \Omega(M(\lambda))$ if and only if
the product
of the matrices $\Psi\Psi_1$ is zero.

\bigskip

%{\bf 3.3 } \ 
The matrix $\Psi\Psi_1$ is diagonal, with diagonal entries
\begin{align*} \lambda p_0r_0(A_0a_0) + (-q_0)s_0(B_0b_0), \cr 
p_1r_1(A_1a_1) + (-q_1)s_1(B_1b_1), \cr
\ldots \cr
p_tr_t(A_ta_t) + (-q_t)s_t(B_tb_t)
\ldots
\end{align*}
for $t\leq m$. Substitute $(-q_t)B_tb_t = p_t(A_ta_t)$ and cancel. 
It follows that:
\bigskip

\begin{lem} We have ${\rm Im}(\Psi_1)={\rm Ker}(\Psi)$ if and only if
$\lambda r_0+s_0=0$ \ and  \  $r_t+s_t=0 \ (1\leq t\leq m)$. 
\end{lem}

We assume this now, and 
we identify the image of $\Psi_1$ with $\Omega(M(\lambda))$.

\bigskip

\subsubsection{The module $\Omega^2(M(\lambda))$.}

Let $P_2=P_0 = \oplus_{t=0}^m e_tA$, and define $\Psi_2: P_2\to P_1$ 
to be the map given by left multiplication with a matrix
of the same form as that of $\Psi$, that
is 
$$\left(\begin{matrix} c_0^+A_0 & d_0^+B_0 && 0 \ldots & 0\cr
0  &   c_1^+A_1 & d_1^+B_1 && 0 \cr
\ldots & & & &\ldots \cr
d_m^+B_m & 0 & \ldots & 0& c_m^+A_m\end{matrix}\right).
$$
Here the $c_i^+$ and the $d_i^+$ are again non-zero scalars. 
We may apply  
Lemma \ref{lem:ident}  again to identify the image of $\Psi_2$.
That is, ${\rm Im}(\Psi_2) = M(\mu)$ where $\mu$ is determined by 
the identities
$$\mu\cdot p_0 d_0^+ + q_0c_0^+=0 \ \ \mbox{ and } 
p_td_t^+ + q_tc_t^+ =0 \ (1\leq t\leq m).
$$
That is
$$ c_0^+ = -\mu(p_0/q_0)d_0^+, \ \mbox{ and }  c_t^+ = -(p_t/q_t)d_t^+.
\leqno{(*)}$$
for $1\leq t\leq m$.
\bigskip

We require that ${\rm Im}(\Psi_2)={\rm Ker}(\Psi_1)$. By comparing dimension,
this is again equivalent with 
$\Psi_1\Psi_2=0$. The matrix $\Psi_1\Psi_2$ is diagonal, 
with diagonal entries
\begin{align*} c_0^+r_0(a_0A_0) + d_m^+s_m(b_mB_m)\cr
c_1^+r_1(a_1A_1) +  d_0^+s_0(b_0B_0) \cr
\ldots \cr
c_l^+r_l(a_lA_l) + d_{l-1}^+s_{l-1}(b_{l-1}B_{l-1})
\end{align*}
(for $1\leq l\leq m$). 

We substitute 
$b_{l-1}B_{l-1} = -(p_l'/q_{l}')a_lA_l$, so we require that
\begin{align*}
c_0^+r_0q_m' -d_m's_mp_0'=0 \cr
c_l^+r_lq_{l-1}' - d_{l-1}^+s_{l-1}p_{l}'=0
\end{align*}
(for $1\leq l\leq m$). 

We know from Lemma 3.13 that $\lambda r_0+s_0=0$ and $r_l+s_l=0$ for $1\leq l\leq m$. We may take $s_t=-1$ for all $t$, 
and then $r_0=\lambda^{-1}$ and $r_l=1$ for $1\leq l\leq m$. 
With this, we get $\Psi_1\Psi_2=0$ if and only if 
$$\lambda^- c_0^+ =   -  d_m^+(p_0'/q_0') \ \mbox{ and} \ c_l^+ = - d_{l-1}^+(p_l'/q_0')
\leqno{(**)}$$
(for $1\leq l\leq m$).

\bigskip

{\bf The proof of Proposition 3.2 for $m\geq 1$.}

We take the product of
all identities in (**), and get
$$\lambda^{-1}\prod_{t=0}^m c_t^+ = (-1)^{m+1}\prod_{t=0}^md_t^+ \cdot
\frac{\prod_{t=0}^m p_t'}{\prod_{t=0}^mq_t'}
\leqno{(***)}$$
We also take the product
over all identities in (*) and get
$$
\prod_{t=0}^m c_t^+ = \mu\cdot (-1)^{m+1}\prod_{t=0}^md_t^+ \cdot
\frac{\prod_{t=0}^m p_t}{\prod_{t=0}^mq_t}
$$
Substitute this into 
(***) and cancel, and we get 
$\mu = \lambda v$
where $v$ is the number in the statement of  Proposition 3.2. 
This proves that
$$\Omega^2(M(\lambda)) \cong M(\lambda v)
$$
$\Box$

\vspace*{1cm}

\section{The proof of the Theorem}

Let $\La_{\bf q}$ be special biserial weakly symmetric, and let $\sigma$
be the permutation of the arrows such that $\alpha\cdot \sigma(\alpha)$ is non-zero in the algebra.
Write $\sigma$ as a product of disjoint cycles. 

\ \ We define the Brauer graph of $\La_q$ as follows. It is the undirected graph whose vertices are the cycles
of $\sigma$. Let $\sigma_1$ and $\sigma_2$ be two cycles of $\sigma$. then the 
edges between $\sigma_1$ and
$\sigma_2$ are labelled by the crossings of $\sigma_1$ and $\sigma_2$. These are the  vertices $i$ of $\cQ$ such that
both $\sigma_1$ and $\sigma_2$ pass through $i$ (counted with multiplicities). 
There is a cyclic ordering of the edges adjacent to a given vertex $\sigma_i$ of the graph;  the successor of edge $e$ is edge $f$ 
if $f$ comes next  after $e$ along the path in $\cQ$ given by $\sigma_i$. This graph is connected, and is independent of ${\bf q}$, we denote it by $G_{\La}$. 

Note that the edges of $G_{\La}$ only see the vertices of $\cQ$ of valency two; we do not need details about vertices
with valency one.  This means that this graph is slightly different from the usual definition of a Brauer graph,  where vertices of $\cQ$ with valency one
are also recorded, the corresponding   edges  $e$ of the graph have the property that one of the adjacent vertices is adjacent only to this edge $e$. 
Hence  our graph is a tree if and only if the usual Brauer graph is a tree. 

Note also that once we know that the graph $G_{\La}$ has no multiple edges and no cycles
then it must be a tree. As well, we do not need to go into details about the cyclic ordering around a vertex.

\bigskip

\begin{expl}
\normalfont 
(1) \ Let $\La_{\bf q}$ be a `Double Nakayama algebra' with
$n$ vertices for $n\geq 2$, where $\sigma$ is a
 product of disjoint 2-cycles. 
That is, $\La_{\bf q} = K\cQ/I$ where $\cQ$ is the quiver
$$
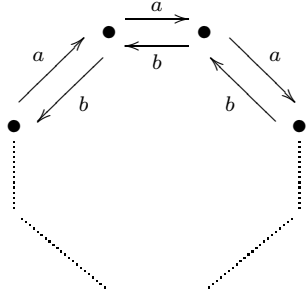
\xymatrix{
& \bullet\ar@<1ex>[r]^a\ar@<1ex>[dl]^{b} & \bullet\ar@<1ex>[l]^{b}\ar@<1ex>[dr]^a & \\
\bullet\ar@<1ex>[ur]^a\ar@{.}[d] & & & \bullet\ar@<1ex>[ul]^{b}\ar@{.}[d]\\
\ar@{.}[rd] & & & \ar@{.}[dl] \\
&&&
}$$
and we label the vertices by $\bZ_n$ and the arrows are $a_i: i\mapsto i+1$ and
$b_i: i+1\mapsto i$. The ideal $I$ is generated by
$a_{i+1}a_{i}, \ \ b_ib_{i+1} $
and
$$p_i(a_ib_i)^{r_i}  + q_i(b_{i-1}a_{i-1})^{r_{i-1}} 
$$
(for $i\in \bZ_n$, where $r_i\geq 1$). 
Then the Brauer graph has $n$ vertices and is a cycle. 
When $r_i=1$ for all $i$ so that  the radical has cube zero,  some socle deformation
does  not satisfy (Fg), by \cite{ES}. By Theorem 1.1
 this holds for arbitrary $r_i\geq 1$.
One can show that for an arbitrary special biserial weakly symmetric
algebra with the above quiver, the Brauer graph is a cycle. 

\medskip

(2) \ Let $\La_{\bf q}$ be an algebra whose quiver is of type $\tilde{Z}$ (with the notation
of \cite{ES}), and where $\sigma$ is a product of 2-cycles together
with 1-cycles for the two loops. That is, $\La_{\bf q} = K\cQ/I$
where $\cQ$ is the quiver 
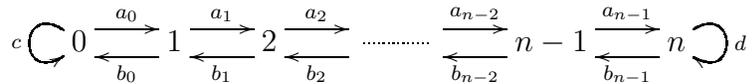
\[
\xymatrix{
0\ar@(ul,dl)_c \ar@<1ex>[r]^-{a_0} &
1\ar@<1ex>[l]^-{b_0}\ar@<1ex>[r]^-{a_1} & 
2\ar@<1ex>[l]^-{b_1}\ar@<1ex>[r]^-{a_2}  
& \ar@<1ex>[l]^-{b_2}\ar@{..}[r] & \ar@<1ex>[r]^--{a_{n-2}}&
n-1\ar@<1ex>[l]^--{b_{n-2}}\ar@<1ex>[r]^--{a_{n-1}} &
n\ar@<1ex>[l]^--{b_{n-1}}\ar@(ur,dr)^d
}\]
where $I$ is generated by the
following relations (we assume $n>0$).
$$
ca_0, \ b_0c, \ a_{n-1}d, \ db_{n-1}, \  a_ia_{i+1}, \ b_ib_{i-1}
$$
$$p_0c^2 + q_0(a_0b_0)^{r_0}, \ \ \ p_i(b_{i-1}a_{i-1})^{r_{i-1}} 
+ q_i(a_ib_i)^{r_i} \  \ \ p_n(b_{n-1}a_{n-1})^{r_{n-1}} + q_nd^2
$$
where $1\leq i\leq n-1$, and $r_i\geq 1$. The coefficients $p_0, \ldots, p_n$
and $q_0, \ldots, q_n$ are non-zero scalars.
Then the Brauer graph of $\La_{\bf q}$ is a line.
When $r_i=1$ for all $i$, 
 it was shown in \cite{ES} that
one can modify the presentation and have  all scalar parameters equal $\pm 1$.
By the Lemma  below,
this holds for arbitrary $r_i$. 

When $r_i=1$ for all $i$, the result of \cite{ES} shows that the algebra
satisfies (Fg). It would be interesting to know whether it is always
the case.

\bigskip

(3) \ Let $\La_{\bf q}$ be the local algebra as in 3.10(1). 
Then $\sigma = (x \ y)$ and hence the Brauer graph has
one vertex with a double edge.
Hence by the Theorem, for some ${\bf q}$, the algebra has criminals.

\bigskip

(4) \ Let $\La_{\bf q}$ be the commutative local algebra as in 3.10(2).
 Then $\sigma  = (x) (y),$ the product of two cycles each of length one, and 
the Brauer graph has two vertices and one edge between them.
Of course we can see here directly that if we rescale generators, then the
scalar parameters in the socle relations
can be changed to $1$ (or anything non-zero). 
\end{expl}

\bigskip

\begin{rem} \normalfont We assume $\La$ has infinite type, then the Brauer graph cannot be just one vertex: If so then  the permutation $\sigma$ would be one cycle with no self-crossings.
 For such an algebra,  all vertices have
valency one and it is  a Nakayama algebra, of finite type.
\end{rem}

\bigskip

We first prove the implication (b) $\Rightarrow$ (a) of Theorem 1.1.
When the Brauer graph of the algebra is a tree,  one can always
rescale the arrows and achieve that all scalar parameters are equal to $1$
(see also the example (2) above). Namely we have the following.

\begin{lem}\label{scale} Assume $\La_{\bf q}$ is a  weakly symmetric and special biserial algebra
whose Brauer graph $G_{\La}$ is a tree. Then $\La_{\bf q}$ is isomorphic to $\La_{\bf 1}$, the algebra where
all parameters are equal to $1$. The algebra $\La_{\bf q}$ does
not have criminals.
\end{lem}

\bigskip

{\it Proof } We will show that by rescaling some arrows one can achieve that all socle parameters become $1$. Note that rescaling arrows does not change the
zero relations of length two.

We fix a cycle $\sigma_0$ of $\sigma$ which has only one neighbour in $G_{\La}$. For any vertex $\sigma_i$ of the Brauer graph, there is
a unique path in $G_{\La}$ of shortest length from $\sigma_0$ to $\sigma_i$. Define the 'distance' 
$d(\sigma_i)$ to be the number of edges of this path. Note also that if $d(\sigma_i)>0$ then $\sigma_i$ has unique neighbour $\sigma_t$ with
$d(\sigma_t) = d(\sigma_i)-1$.

We prove the Lemma by induction on the distance. If $d=0$ then the cycle is $\sigma_0$, and we keep its arrows as they are.

For the inductive hypothesis, assume that for all cycles with $d(\sigma_t)< d$, the arrows in it have been scaled so that the relevant socle relations have
scalars equal to $1$.  Now take $\sigma_i$ such that $d(\sigma_i)=d+1$. Then there is a unique $\sigma_t$ joined in the Brauer graph to $\sigma_t$
such that $d(\sigma_t)=d$. Let
$j$ be the edge between $\sigma_t$ and $\sigma_i$ in the Brauer graph. That is, $j$ is a vertex in $\cQ$ of valency two.
Consider the socle relation at $j$,
$$p_jC_j + q_jD_j=0$$
Say the arrows in $C_j$ are the arrows of $\sigma_t$, so that $C_j = (\alpha_1\alpha_2\ldots \alpha_r)^{m_t}$ where
$\sigma_t =   (\alpha_1 \ \alpha_2 \ \ldots \alpha_r)$ and $m_t\geq 1$.
Then $D_j = (\beta_1\beta_2\ldots \beta_s)^{m_i}$ where $\sigma_i = (\beta_1 \ \beta_2 \ \ldots \beta_s)$ and $m_i\geq 1$.
We replace a single arrow in $\sigma_i$, namely $\beta_1$, by $\beta_1':= c\beta_1$ where $c$ is a root of  $x^{m_i} - (q_i/p_i)$. Then $\beta_1'$ is an arrow, and replacing
$\beta_1$ by $\beta_1'$ does not affect zero relations of length two.  The monomial $D_j':= (\beta_1'\beta_2\ldots \beta_s)^{m_i}$ is an element in the socle.
By the choice of $c$ it follows that
$$p_j(C_j + D_j') = 0 \ \ \mbox{ and hence} \ C_j + D_j'=0.
$$
Since we have not changed any of the $\alpha_u$'s, the relations fixed earlier are not
altered.

\ \ We may take all parameters
equal to $1$, and then by Proposition 
3.7, the module $M(\lambda)$ for any band word is periodic.
By the discussion in Section 2.1, the algebra
$\La_{\bf 1}$ does not have any criminals at all.
$\Box$

\bigskip

We observe that we could equally well have signs $1$ and $-1$. With this
one can show that the algebra is symmetric if its  Brauer graph is a tree.

\bigskip

For the implication (a) $\Rightarrow$ (b) of Theorem 1.1, we start with the following.
We will use  band words, as defined in 3.1, which involve
vertices of valency two. When the word $W$ is fixed, we 
write $\cE=\{ e_0, \ldots, e_m\}$ and
$\cF=\{ f_0, \ldots, f_m\}$ for the sets of these vertices (which depend on
$W$). 

\bigskip

\begin{lem}\label{e-not-f} Assume $\La_{\bf q}$ has a band word $W$ for which 
$\cE \neq \cF$ and where 
the $f_i$ are distinct and the $e_i$ 
are distinct. Then for some choice of ${\bf q}$ the algebra has
criminals.
\end{lem}

\bigskip

{\it Proof } Suppose, say,  $f_i\not\in \cE$. The socle relation 
at $f_i$ which we denoted by $(\theta_i)$ contributes the factor
$q_i/p_i$ to the parameter $v$ of the word $W$. 
The relation $(\theta_i)$ does not occur elsewhere since the
$f_j$ are distinct and $f_i$ is not in $\cE$. Take $q\in K$ which is not
a root of unity, and take $q_i:=q$ and set all other parameters for
socle relations equal to $1$. Then $v=q$ and hence a module $M(\lambda)$ with
this word is a criminal fo $\La_{\bf q}$.

\bigskip

\begin{prop}\label{e=f} Assume we have a band word \ $W$ where 
$\cE = \cF$, of size $m+1$, and let $\pi$ be the permutation of
$m+1$ with $e_{\pi(i)} = f_i$. Then the following
are equivalent:\\
(1) \ The permutation $\sigma$ takes the last arrow of $a_i$ to the 
first arrow of $a_{\pi(i)}$ for all $i$;\\
(2) \ $v=1$  where $v$ is the parameter for $W$.
\end{prop}

{\it Proof } \ Fix a vertex $i$, the socle 
relations for $f_i$ and for $e_{\pi(i)}$ are equal.
Recall we have written them as 
$$p_i(A_ia_i) + q_i(B_ib_i)=0,  \ \ \
 p_{\pi(i)}'(a_{\pi(i)}A_{\pi(i)}) + q_{\pi(i)}'(b_{\pi(i)-1}B_{\pi(i)-1}) = 0.
$$
For the moment, view $p_i, q_i$ and $q_{\pi(i)}', p_{\pi(i)}'$ as
indeterminates.
The contribution of these to the parameter $v$ is 
$$\frac{q_i}{p_i}\cdot \frac{p_{\pi(i)}'}{q_{\pi(i)}'}.
$$
Hence $v=1$ if and only if for all $i$ we have
$q_i=q_{\pi(i)}'$ and
$p_i=p_{\pi(i)}'$. 
This holds if and only if for all $i$ we have
$A_ia_i = a_{\pi(i)}A_{\pi(i)}$, (or equivalently
$B_ib_i = b_{\pi(i)-1}B_{\pi(i)-1}$). 

\medskip

Recall that a rotation of the path $A_ia_i$ is a non-zero element in the
algebra. Therefore if $\alpha$ is the last arrow of $a_i$ then 
$\sigma(\alpha)$ must be the first arrow of $A_i$.

Since $A_ia_i$ is either $a_{\pi(i)}A_{\pi(i)}$ or $b_{\pi(i)-1}B_{\pi(i)-1}$, the
arrow $\sigma(\alpha)$ is the first arrow of precisely one of $a_{\pi(i)}$ or
$b_{\pi(i)-1}$.

Hence $A_ia_i = a_{\pi(i)}A_{\pi(i)}$  if and only if 
$\sigma(\alpha)$ is the first arrow of $a_{\pi(i)}$.

$\Box$

\bigskip

{\bf  Proof of (a) $\Rightarrow$ (b) of  Theorem 1.1. } \ 
 
 Assume that for any ${\bf q}$ the algebra $\La_{\bf q}$ does not
have criminals. We must show that the Brauer graph $G_{\La}$ is a tree
with no multiple edges.

Take a band word $W$ as in 3.1 in which all paths $a_i$ and $b_i$ have
minimal length, that is, all vertices along these paths other than $e_i$ and
$f_i$ (if any) have valency one. Such $W$ must exist, we refer to the $a_i$
and $b_i$ as minimal paths in this proof.

(1) \ We claim that $\cE = \cF$. Note that by the minimality condition
in 3.1 the vertices 
$e_0, \ldots, e_m$ are pairwise distinct, similarly the $f_0, \ldots, f_m$
are pairwise distinct.
Hence the claim follows from  Lemma \ref{e-not-f}.

\bigskip
(2) \ We claim that  $\cE$ is the set of all vertices of valency two  of $\cQ$:

Suppose not. Since $\cQ$ is connected, there must be a minimal path $\gamma$ say starting or ending at some vertex of valency two $e$ which is not in $\cE$, and
ending or starting at some vertex $e_i\in \cE$. But there are only two 
minimal paths  starting at $e_i$ and two minimal paths ending at $e_i$, and $\gamma$ must then be one of the $a_i$ or $b_i$, a contradiction.

\bigskip

(3) \ We claim that the set of arrows occuring in the paths $\{a_0, a_1, \ldots, a_m\}$ is invariant under the permutation $\sigma$: 

If $\alpha$ is an arrow ending at a vertex $j$ of valency one then clearly $\sigma(\alpha)$ is 
the arrow starting at $j$. So we only need to know that $\sigma$ takes the last arrow
of some $a_i$ to the the first arrow of some $a_l$. 
We have no criminals, therefore $v=1$ and the claim holds by Proposition 4.5.

Then as well, $\sigma$ leaves  the set of arrows invariant 
occuring in any of the paths $\{ b_0, b_1, \ldots, b_m\}$.

\bigskip

This means that
we can colour the cycles of $\sigma$ by the two colours $a$
and $b$, and this gives a colouring for the vertices  of  $G_{\La}$.

\bigskip

(4) \ By considering the word $W$ we see that  that for each vertex $i$ of  valency two, there is
one $a$-cycle and one $b$-cycle passing through $i$. 
Hence the Brauer graph does not have edges between two cycles
of the same colour,  and there is no edge in $G_{\La}$ 
starting and ending at the same
cycle of $\sigma$.

\medskip

(5) \ The graph $G_{\La}$ does not have multiple edges: 
Suppose $\sigma_1, \sigma_2$ are cycles which pass through 
vertices $e \neq f$. 
Then we can find a band word $W$ of the form
$$e \stackrel{a_0}\rightarrow f 
\stackrel{b_0}\leftarrow e$$
Namely, say $\sigma_1$ is an $a$-cycle, then take
for $a_0$ the shortest path consisting of arrows in $\sigma_1$ from $e$ to $f$, and take $b_0$ similarly. Then the parameter $v$ for this word is
$\neq 1$, by the Lemma \label{e-not-f}.

So far, we have proved that the Brauer graph has a colouring of vertices with colours $a$ and $b$ where the colours alternate, and 
it does not have multiple edges.

(6) \ The graph $G_{\La}$   does not have a cycle:  

If there is a cycle in the graph,  then this cycle
 must have an {\it even } number of vertices, since
the vertices of $G_{\La}$ are coloured by two alternating colours.
If we take part of each cycle in the appropriate direction, then we get a band word with $\cE$ and $\cF$  disjoint,
since the number of vertices is even. 
By  Lemma 4.4, this gives rise to a criminal.
We have now proved that $G_{\La}$ is a connected graph with no multiple edges, and without cycles, and hence $G_{\La}$ must be a tree.
$\Box$

\medskip

{\sc Karin Erdmann, Mathematical Institute, University of Oxford,
ROQ, Oxford OX2 6GG, UK}

\textit{email:} {\tt erdmann@maths.ox.ac.uk}

\end{document}